\documentclass[smallextended, numbook]{svjour3}       
\smartqed  
\input{microtyping}

\usepackage[english]{babel}
\usepackage{booktabs}
\usepackage{multirow}
\usepackage{color}
\usepackage[table]{xcolor}
\usepackage{tikz}
\usepackage{pgfplots}
\usetikzlibrary{plotmarks}
\newcommand\blue[1]{{\color{black}#1}}
\newcommand\green[1]{{\color{black}#1}}
\newcommand\greenNew[1]{{\color{black}#1}}

\usepackage{placeins}
\usepackage{siunitx}
\usetikzlibrary{arrows.meta}
\usepackage{graphicx}
\usepackage{multicol}
\usepackage{algorithm}
\usepackage[noend]{algpseudocode}
\usepackage[authoryear]{natbib}
\usepackage{amsfonts,amsmath,amssymb}
\allowdisplaybreaks 
\usepackage{mathtools}
\newcommand{\csum}[1]{\sum_{\mathclap{#1}}}
\newcommand{\set}[2]{\left\{\, #1 \mathrel{}\middle|\mathrel{} #2 \,\right\}}
\DeclareMathOperator*{\argmax}{argmax}
\DeclarePairedDelimiter\floor{\lfloor}{\rfloor}
\DeclarePairedDelimiter\ceil{\lceil}{\rceil}

\usepackage[utf8]{inputenc}
\usepackage{graphicx}
\usepackage{verbatim}
\usepackage[autostyle, english=american]{csquotes}
\MakeOuterQuote{"}

\usepackage{url}
\usepackage{subfig}
\usepackage{dblfloatfix}
\usepackage{cite}
\usepackage{setspace}
\usepackage{multicol}

\pgfplotsset{compat=1.15}

\newcommand\eg{e.\,g.}
\newcommand\ie{i.\,e.}
\newcommand\brandec{\textsc{BranDec}\footnote{\url{https://gitlab.com/Soha/brandec}}}
\newcommand\ver{v0.8.4}
\newcommand{\named}[6][\min]%
{ \begin{subequations}%
\label{mod:#2}%
\begin{align}%
#1   && #4 \span  \tag{\textbf{#3}} \label{tag:#2} \\
s.t. && #5 \\
&& #6 \span  \nonumber
\end{align}%
\end{subequations}}
\renewcommand\exp[1]{\times 10^{#1}}

\usepackage[%
rm={oldstyle=false,proportional=true},%
sf={oldstyle=false,proportional=true},%
tt={oldstyle=false,proportional=true,variable=false},%
qt=false%
]{cfr-lm}
\usepackage[noabbrev, capitalise]{cleveref}

\begin{document}

\title{Solving the Integrated Bin Allocation and Collection Routing Problem for Municipal Solid Waste: a Benders Decomposition Approach
}

\titlerunning{A Benders approach for location-routing in tactical MSW}        

\author{Arthur Mahéo         \and
        Diego Gabriel Rossit \and
        Philip Kilby
}

\authorrunning{Mahéo et al.} 

\institute{Arthur Mahéo \at
              Monash University, Faculty of IT \\
              Melbourne, Australia \\
              \email{arthur.maheo@monash.edu}           
           \and
           Diego Gabriel Rossit \at
              Department of Engineering, Universidad Nacional del Sur \\
              Instituto de Matem\'{a}tica de Bah\'{i}a Blanca, CONICET \\
              \email{diego.rossit@uns.edu.ar}  
            \and
            Philip Kilby \at
              Australian National University \\
              \email{Philip.Kilby@data61.csiro.au}  
}

\date{Received: date / Accepted: date}

\maketitle

\begin{abstract}
The municipal solid waste system is a complex reverse logistic chain which comprises several optimisation problems.
Although these problems are interdependent -- i.e., the solution to one of the problems restricts the solution to the other -- they are usually solved sequentially in the related literature because each is usually a computationally complex problem.
\green{We address two of the tactical planning problems in this chain by means of a Benders decomposition approach: determining the location and/or capacity of garbage accumulation points, and the design and schedule of collection routes for vehicles.}
Our approach manages to solve medium-sized real-world instances in the city of Bahía Blanca, Argentina, showing smaller computing times than solving a full MIP model.
\end{abstract}

\section{Introduction}
\label{intro}

\blue{Regardless of their size, city councils have the duty to provide efficient service to their constituents. Municipal Solid Waste (MSW) management is one such a crucial service. When mishandled, this problem may produce serious economic, environmental and social impacts~\citep{asefi2019mathematical,hoornweg2012waste}.}

In this paper, we will focus on a less traditional MSW design called Garbage Accumulation Points (GAPs).
Instead of providing a ``door-to-door'' pickup of garbage, constituents have to drop their garbage at specific facility -- the GAPs. These facilities can range from collective community bins to recycling centres.

\green{Although the decision between using a ``door-to-door'' or a ``GAP-based'' waste collection system is site-specific~\citep{rossit2022wastebin}, there are recent studies that showed that the GAP-based system is more efficient in terms of transporting cost since the distances travelled by vehicles are reduced~\citep{blazquez2020network} which can also lead to a smaller environmental impact through the reduction of  greenhouse gas (GHG) emissions and other air pollution-related metrics~\citep{gilardinoX2017combining}. The savings in the overall MSW system that can be achieved with a reduction in the transportation costs can be even more important in countries that experience relatively high logistic costs, such as Argentina~\citep{brozX2018argentinian,alalog2021}, where the computational experimentation of this work is performed. }

When using GAPs, MSW management comprises the following design decisions:
\begin{itemize}
\item The design of a pre-collection network, which consists in defining the location and capacity of GAPs.
\item The design and schedule of routes for collection vehicles.
\end{itemize}
The geographical distribution of GAPs affects the actual route that the collection vehicles must perform.
Additionally, the storage capacity of these sites will define the visit
frequency in order to avoid overflow.
Finally, the availability and type of vehicles\footnote{Mainly capacity, but could also be cost.} affect the distribution and capacity of the GAPs in a (global) optimal solution.
Thus, there is a trade-off between the cost of the installation of GAPs and the routing cost; solving both simultaneously is often beneficial~\citep{hemmelmayrX2013models}.

However, solutions in the literature \citep{ghianiX2014operations,rossit2022wastebin,hanX2015waste} often address each separately -- see \cref{sec:review}.
This is due to the complexity of tackling MSW as a whole.
Indeed, only solving the design of routes is tantamount to solving a Vehicle Routing Problem (VRP) \citep{toth2002vehicle}, a well-known NP-hard problem.

In this paper, we propose the following contributions to the field of MSW management with GAPs:
\begin{itemize}
    \item A novel mathematical model which combines the allocation of bin combinations to GAPs and defining collection routes (\cref{sec:mathematical_model}).
    \item A Benders decomposition-based approach to tackle the resulting problem (\cref{sec:benders}).
\end{itemize}
The problem we are tackling is an ``inventory routing problem'' \citep{campbell1998inventory}.
The resulting formulation is a mixed-integer program (MIP) which is still too large to be tractable.
However, we can see it as a combination of two problems:
\begin{enumerate}
    \item a routing problem, similar to a vehicle routing problem \citep{toth2002vehicle}; and,
    \item an allocation problem, similar to a nonlinear resource allocation problem \citep{bretthauer1995nonlinear}, in which the used amount of resource (bin) should be minimised, though not limited.
\end{enumerate}
This \emph{natural} decomposition lead us to use Benders decomposition \citep{benders1962partitioning}, a well-suited method for problems with this structure.
Benders decomposition works by solving such problems in an iterative fashion.
First, it solves the \emph{difficult} part to generate a candidate solution.
It then checks this solution against the dual of the \emph{easy} part.
From the dual solution, it either terminates, when the dual solution's objective value is equal to an incumbent; or, it generates constraints, called ``Benders cuts,'' which are added to the difficult part and the problem is solved anew.

However, we cannot use standard Benders decomposition because the subproblem contains integer variables.
Therefore, we use a framework called Unified branch-and-Benders-cut  \citep[UB\&BC,][]{maheo2020unified}.
This framework is based on a modified Branch-and-Cut (B\&C) with callbacks from a commercial solver.
In the callbacks, it derives dual information and an upper bound for the subproblem.
Using these, it terminates the branch-and-bound tree with a set of \emph{open solutions} -- whose objective function value falls below the best upper bound.
To find the global optimum, the B\&C is followed by a post-processing phase where the framework solves those open solutions to integer optimality.

We test our model on a real-world use case: the city of Bahía Blanca, Argentina (\cref{sec:experimentation}).
Although the city currently uses a door-to-door collection service, they are interested in switching to GAPs. We simulated instances using data from a survey \citep{cavallinX2020application} and provide optimal allocation and routing for a variety of scenarios. \blue{Preliminary results of this work were presented at the  \textit{10\textsuperscript{th} International Conference of Production Research -- ICPR Americas 2020}~\citep{maheo2020benders}. The new content in this article include a preprocessing algorithm for discarding inefficient bin combinations in advance, an improved set of realistic instances with a larger documentation of how relevant data is gathered, an updated and comprehensive literature review and new features to enhance the Bender's resolution process. Additionally, several sections were completely rewritten to enhance readability of the manuscript.}

\green{
This work is structured as follows. In \cref{sec:review} we present an updated review of the related works in the literature. In \cref{sec:mathematical_model} we present a mathematical formulation of the problem improved with valid inequalities. In \cref{sec:benders} we present the resolution approach based on Benders decomposition. In \cref{sec:example} we proposed an illustrative working example to clearly outlined our resolution algorithm. In \cref{sec:experimentation} we present the computational experimentation. Finally, we present our conclusions and future directions in \cref{sec:conc}.}

\section{Literature review}\label{sec:review}

Allocation of bins and routing problems have been thoroughly studied as separate problems in the MSW related literature~\citep{luX2015smart}.
Comprehensive reviews of the study of these problems separately can be found in \citet{rossit2022wastebin} for the allocation of bins problem and in \greenNew{\citet{banyai2019optimization}} and \citet{hanX2015waste} for the routing problem.
However, the number of works considering integrated approaches is more scarce. In this Section, we present the main related works according to three different categories: works that are related to integrated approaches to collect unsorted waste, works that are related to integrated approaches to collect recyclable material, and works that used Benders' decomposition in other stages of the MSW reverse logistic chain.

\green{
\subsection{Integrated approaches to collect unsorted waste}

\Citet{hemmelmayrX2013models} proposed an integrated approach where the bins allocation problem is solved jointly with the routing schedule.
The authors compare a Variable Neighbourhood Search (VNS) algorithm for solving the problem hierarchically -- \ie, first solving the bin allocation and then the routing and \emph{vice versa} -- and integrated approaches.
They found that integrated approaches give better results than hierarchical ones.
The same strategy  -- \ie, comparing integrated approaches with hierarchical approaches -- was implemented in \citet{kim2015integrated} for a locating routing problem solved with a Tabu Search algorithm.
Computational experimentation on benchmark instances and a real case study of Seoul, South Korea~\citep{kim2015case}, found that integrated approaches allowed obtaining better solutions.
Another example is \citet{jammeliX2019bi} who presented a study case of the Tunisian city of Sousse, considering uncertainty in waste generation at GAPs.
They consider that all GAPs are to be collected daily.
They proposed a transformed formulation to handle stochastic waste generation and solved the problem in a heuristic fashion: first, they applied the $k$-means clustering algorithm to group the GAPs into sectors and, later, they applied an exact model solved with CPLEX to determine both the number of bins and the collection route of each sector.

\subsection{Integrated approaches to collect recyclable material}

Another popular area for location-routing applications is the collection of recyclable materials.
\Citet{chang1999strategic} propose an integrated approach.
They used an evolutionary algorithm to solve a location-routing problem for the city of Kaohsiung, Taiwan.
Their approach is multi-objective as they maximise the population serviced, and minimise the total walking distance, from household to recycling drop-off stations, and the total driving distance of the collection vehicles.
This work is extended in~\citet{chang2000siting}, where the same three objective functions are considered as fuzzy goals.
\Citet{vidovic2016two} presented an integrated approach for the location-routing problem proposing an integrated MIP formulation and a hierarchical two-step heuristic approach.
Similarly to~\citet{yaakoubi2018heuristic}, the hierarchical approach is able to obtain near optimal solutions, comparable to the integrated MIP formulation, for synthetic instances designed by the authors.
Another similar case was proposed by~\citet{sheriff2017integrated} in which they compared an integrated MIP formulation with a multi-echelon heuristic.
However, as opposed to the work of~\citet{vidovic2016two}, the integrated approach obtained better results than the sequential approach on computational tests performed in an (non-specified) Indian urban area.
Another approach is presented in \citet{hemmelmayrX2017periodic} for solving an integrated model that aims to simultaneously locate GAPs, size the storage capacity of each GAP (allocate bins) and set the weekly collection schedule and routes in the context of collaborative recycling problem.
They solved this problem with an Adaptive Large Neighbourhood Search algorithm based on their previous implementation \citep{hemmelmayr2015sequential}.
They performed a sensitivity analysis for several of the parameters, such as available vehicle capacities, visiting schedules or GAP storage capacities.
Finally, \citet{cubillos2020solution} presented a location-routing approach for recyclable material where they approximate the collection with a Travelling Salesman Problem -- \ie, they considered an uncapacitated vehicle that visits all the bins.
For solving the integrated problem, they proposed a Variable Neighbourhood Search (VNS) algorithm.
They used that method to solve real instances of four Danish cities.
\Citet{gultekin2020decomposition} present a location-routing application for a specific type of waste: cooking oil that is collected to produce biodiesel.
They proposed an integrated MIP to simultaneously define the location of GAPs, the assignment of waste generators to the available GAPs and define the collection routes.
They present a heuristic to address this problem dividing it in two parts.
The first part solves the GAP location and generators assignment.
The first part generates different routes that are later used in the second part which solves the routing schedule.
The problem remains too difficult to solve to optimality, so they save solutions found until a time limit.
This approach is competitive when compared to their MIP model in a set of synthetic instances.

\subsection{Benders' decomposition in other stages of the MSW reverse logistic chain}

Other approaches using Benders decomposition exist to deal with optimisation problems of the strategic level of the MSW logistic chain.
Most consider stochastic parameters, for which Benders decomposition is a traditional application area.
For example, \citet{saif2019municipality} used Benders decomposition to model a logistic chain of MSW in which organic waste is sent from sources to treatment plants to generate power.
They consider Uncertainty in waste generation, power price, and demand.
Another case is presented in~\citet{kuudela2019multi}, who applied Benders decomposition to optimise the location and capacity selection of waste transfer stations when considering uncertainty in the operational cost of the stations.
\citet{Fattahi2020data} proposed a data-driven stochastic programming model based on Benders Decomposition that is applied to a case study in Tehran, Iran.
The aim is to design a MSW recovery network for power generation considering uncertain waste generation rate.

\subsection{Summary and contribution of this work}

After having revised the main related works we acknowledge that there is still room to propose efficient methods for solving this complex locating-routing problem in the initial stages of the MSW reverse logistic chain. Moreover, the used of Benders' decomposition, which have been efficiently applied in other complex MIP models, has not been previously used to address this particular problem.
}

\section{A mathematical model of MSW}
\label{sec:mathematical_model}

In this section we present a mixed-integer mathematical formulation for the integrated problem of simultaneously bin combinations to GAPs and defining collection routes in the context of the MSW system.
\blue{The formulation takes advantage of a preprocessing phase that establishes the use of bins combinations.}
We improve the formulation with the addition of a set of valid inequalities.

\blue{%
\subsection{Bin combinations preprocessing}
\label{sec:preproc}

Realistic problems usually involve locating different types of bins.
The bins have different purchasing and maintenance cost, storage capacity (for accumulating waste) and occupied space.
Usually a GAP has enough space to locate more than one bin -- \ie, a bin combination.
There are two main aspects to consider when deciding which bin combination can be installed in a given GAP:
\begin{description}
  \item[A feasibility constraint:] GAPs have a maximum available space and, thus, the bin combination has to fit in that limited space~\citep{toutouh2020soft}.
  \item[An efficiency criterion:] we can avoid economically inconvenient bin combinations.
\end{description}

The last point comes from the fact that bin combinations will have different characteristics.
First, bin combinations will have different \emph{joint storage capacity} -- \ie, the sum of the capacities of the bins that conforms the bin combination.
They will also have different \emph{joint cost} -- \ie, the sum of the purchasing and maintenance costs.
Thus, some combinations will have a larger joint cost and a smaller joint capacity than other.
These bin combinations can be dismissed since they will not appear in any \emph{efficient} solution.
We can therefore determine the set of feasible and Pareto-optimal bin combinations for each GAP. 
We developed a preprocessing algorithm (\cref{Alg:EA}) to allow a more compact mathematical formulation of the problem.

\begin{algorithm}[!ht]
\small
\caption{Preprocess Pareto-optimal bin combinations}
\label{Alg:EA}
\begin{algorithmic} [1]
\Function{\textsc{Preprocess}$(P, sp_{GAP})$}{}
    \State {\bf Initialise} list $L$ as empty
    \State $s_{min} = \min_{i \in P} s(i)$ 
\State $r' = \floor*{\frac{sp_{GAP}}{s_{min}}} $
    \For{$r = 1; r \leq r'; r++$}
        \State {\bf Initialise} container $V$ as empty and size $r$
        \State {\bf Call} {\textsc{Fill}($V,P,0,r,0,|P|-1,L,sp_{GAP}$)}
    \EndFor
    \State Return $L$
\EndFunction

\State 
\Function{\textsc{Fill}$(V,P,index,r,start,end,L,sp_{GAP})$}{}
    \If{$index == r$}
        \If{$s(V) \leq sp_{GAP}$} \Comment{function $s()$ measures the occupied space}
        \State \textsc{Update}$(L,V)$
        \EndIf
    \EndIf
    \For{$i = start$; $i \leq end$;$i++$}
        \State {$V[index] \leftarrow i$}
        \State \textsc{Fill}$(V,P,index+1,r,start,end,L,sp_{gap})$
    \EndFor
\EndFunction

\greenNew{
\State 
\Function{\textsc{Update}$(L,V)$}{}
\State {Calculate joint storage capacity and joint cost of $V$}
\State {$H =$ sorted($L \cup V$, descending order joint storage capacity)}
\State {$L =$\textsc{Front}$(H)$}
\EndFunction

\State
\Function{\textsc{Front}$(H)$}{}
    \If{$|H| == 1$}
        \State{Return $H$}
    \Else
    \State{$T =$ Front$(H[1:|H|/2])$}         \State{$B =$ Front$(H[|H|/2 + 1: |H|])$}
     \State {\bf Initialise} list $M$ as empty
    \For{$b$ in $B$}
        \If{$b$ is not dominated by any bin combination in $T$ regarding joint cost}
        \State{$M = M \cup b$}
        \EndIf
    \EndFor
    \State{Return $M \cup T$}
    \EndIf
\EndFunction
}

\State 
\end{algorithmic}
\end{algorithm}

The \textsc{Preprocess} receives as an input the set of types of bin ($P$) and the available space in the GAP ($sp_{gap}$) and returns a list of feasible and Pareto-optimal bin combinations ($L$). In short, the algorithm works as follows:
\begin{enumerate}
    \item It estimates the value $r'$ which is the maximum number of bins that can be placed in the GAP (this is done considering the smallest type of bin $s_{min}$ and the available space of the GAP $sp_{gap}$).
    \item It evaluates the convenience of each possible bin combination $V$ that has a number of bins smaller or equal to $r'$ in two steps:
    \begin{enumerate}
     \item It analyses if bin combination $V$ is feasible considering the available space in the GAP.
     \item If bin combination $V$ is feasible, it applies function \textsc{Update}$()$ to update list $L$ to store only Pareto-optimal bin combinations of the set $L \cup V$ (this is considering the relation between the joint cost and the joint storage capacity of the bin combination).
     \greenNew{This function sorts the set of solutions set $L \cup V$ in terms of joint storage capacity and then uses the well known \citeauthor{kung1975finding}'s method to get the Pareto front of a set of solutions of a multi-objective problem. The detailed outline of this algorithm can be consulted in~\citet{kung1975finding}. 
     In the case of our preprocessing algorithm the Kung's method works as follows. A recursive function \textsc{Front}$()$ is applied to first split the input set of solutions in two halves and then compares the second half with elements of the first half in order to get the non-dominated solutions.}
    \end{enumerate}
    \item It returns list $L$ which stores only the Pareto-optimal bin combinations for the GAP.
\end{enumerate}

Hereafter, we understand that all the bin combinations were computed following this procedure and, thus, are Pareto-optimal.
}

\subsection{Model formulation}
\label{subsec:formulation}

The mathematical model has the following sets:

\begin{itemize}
\item $I$: the set of potential GAPs.
\item $L=\{l_0, l_1, \dots, l_{|L|}\}$: is an ordered set of vehicles. We consider a homogeneous and finite fleet of vehicles.
\item $T$: the set of days in the time horizon, which coincides with a week (seven days).
\item $R$: the set of possible visit combinations.
\item $U$: the set of all bin combinations that can be installed in a GAP.
\end{itemize}

A potential GAP $i\in{I}$ is a predefined location in an urban area in which bins can be installed.
We define the superset: $I^0 = I \cup {0}$, where $0$ is the depot from which vehicles start and finish their daily tours, and where the collected waste is deposited.
We also define a special notation for the set of edges given a set of nodes: 
$$ E(I) = \set{(i, j)}{i \in I, j \in I, i \neq j}$$.
The set of visit combinations $R$ are possible weekly schedules to empty GAPs.
The set of bin combinations $U$ is obtained using \cref{Alg:EA} and represent feasible and Pareto-optimal arrangements.

We now define the parameters of the model:

\begin{itemize}
\item $Q$: vehicle capacity.
\item $c_{ij}$: travel time between $i$ to $j$.
\item $s_i$: service time of a GAP $i$.
\item $b_i$: waste generation per day at GAP $i$.
\item $cap_u$: capacity of bin combination $u$.
\item $cin_u$: adjusted cost of installing bin combination $u$ for the time horizon $T$.
\item $\alpha$: cost per kilometre of transportation.
\item $\beta_r$: maximum number of days between two consecutive visits of the visit combination $r$.
\item $a_{rt}$: 1 if day $t$ is included in visit combination $r$.
\item $TL$: time limit of the working day.
\end{itemize}

Notice that $cin_u$ is an \emph{adjusted cost}.
This is because we are considering two different level of decision and cost:
\begin{enumerate}
\item a strategic decision that involves purchasing and installing the bin combinations that will last probably for several years; and,
\item a tactical decision which involves the transport costs of the routing schedule~\citep{nagy2007location}.
\end{enumerate}
Therefore, the cost assigned to a bin combination ($cin_u$) includes a proportional part of the purchase and installation costs, and the maintenance cost.

With regards to parameters $a_{rt}$ and $\beta_r$, let us introduce them with an example:
\begin{example}
Let the time horizon be a week: $T = \{t_1,t_2,t_3,t_4,t_5,t_6,t_7\}$.
Then, one possible visit combination $r* \in{R}$ is $\{t_1,t_3,t_5,t_7\}$.
In this case, we have: $a_{r* t_1} = a_{r* t_3} = a_{r* t_5} = a_{r* t_7} = 1$, and, conversely: $a_{r* t_2} = a_{r* t_4} = a_{r* t_6} = 0$.
Thus, the maximum number of days between two consecutive visits that this combination has is two days: $\beta_{r*} = 2$, and the chosen bin combination for this GAP must be able to store the waste generated in two days.%
\footnote{A similar consideration is performed in \citet{hemmelmayrX2013models}.}
\end{example}

Finally, we define the following decision variables:

\begin{itemize}
\item $x_{ijlt}$: binary variable set to 1 if vehicle $l$ performs the collection route between GAPs $i$ and $j$ on day $t$, 0 otherwise.
\item $v_{ijlt}$: continuous variable representing the load of vehicle $l$ along the path between GAP $i$ and $j$ on day $t$.
\item $m_{ir}$: binary variable set to 1 if visit combination $r$ is assigned to GAP $i$, 0 otherwise.
\item $n_{ui}$: binary variable set to 1 if bin combination $u\in{U}$ is used for GAP $i$, 0 otherwise.
\end{itemize}

We now present the mathematical model for the MSW management problem:%
\named{m1}{M1}
      {\sum_{\substack{i\in{I} \\ u\in{U}}} n_{ui}\ {cin}_{u} + \alpha \csum{i, j \in E(I^0)} \left(c_{ij} + s_i\right) \left(\sum_{\substack{{l\in{L}} \\ t\in{T}}} x_{ijlt} \right)\span}
      {
\label{eq:2}
\sum_{u\in{U}} n_{ui}\ cap_u
 & \geq \sum_{r\in{R}} b_{i} m_{ir} \beta_r  \span \forall{}\ i\in{I} \\
&& \label{eq:2b}
\sum_{u\in{U}} n_{ui} & =  1& \forall{}\ i\in{I} \\
&& \label{eq:3}
\sum_{r\in{R}} m_{ir} & = 1& \forall{}\ i\in{I} \\
&& \label{eq:4}
\csum{\substack{j \in{I^0},g\neq{i} \\ l\in{L}}} x_{ijlt} - \sum_{r\in{R}} a_{rt} m_{ir} & = 0& \forall{}\ t\in{T}, \ i\in{I} \\
&& \label{eq:5}
\csum{i\in{I^0},i\neq{q}} x_{ijlt} - \csum{j\in{I^0},j\neq{q}} x_{qjlt} &= 0& \forall{}\ q\in{I^0},\ l\in{L},\ t\in{T} \\
&& \label{eq:6}
\sum_{i\in{I}} x_{0ilt} & \leq 1& \forall{}\ l\in{L},\ t\in{T} \\
&& \label{eq:7}
\csum{(i,g) \in E(I^0)} (c_{ij} + s_i)\ x_{ijlt} & \leq TL& \forall{}\ l\in{L},\ t\in{T} \\
&& \label{eq:8}
v_{ijlt} & \leq Q\ x_{ijlt} & \forall{}\ (i, j) \in E(I^0), l\in{L},\ t\in{T} \\
&&
\csum{i\in{I^0}, i\neq{g}} v_{ijlt} + b_g\
\sum_{r\in{R}} \left(m_{gr}\ \beta_r\right)
\leq
\csum{i\in{I^0}, i\neq{g}} v_{jilt} + Q\
\left(1 - \csum{i\in{I^0}, i\neq{g}} x_{ijlt}\right)
\span \span \nonumber \\
&&&&\forall{}\ j \in{I},\ l\in{L},\ t\in{T} \label{eq:9}
      }
      {v \geq 0; m, n, x, b \in \mathbb{B}}

The objective function is the sum of the routing cost and the adjusted cost of installing bins.
\Cref{eq:2} establishes that, given a visit combination, the maximum amount of garbage that can be accumulated in a GAP cannot surpass the installed capacity of the bin combination.
\Cref{eq:2b} enforces that one bin combination has to be chosen for each GAP.
\Cref{eq:3} establishes that one visit combination is assigned to each GAP.
\Cref{eq:4} ensures that each GAP is visited by the collection vehicle the days that corresponds to the assigned visit combination.
\Cref{eq:5} ensures that if a vehicle visits a GAP, it leaves the GAP on the same day.
\Cref{eq:6} states that every vehicle can be used at most once a day.
\Cref{eq:7} guarantees that a tour does not last longer than the allowable time limit associated with the working day of the drivers.
\Cref{eq:8} limits the total amount of waste collected in a tour to the vehicle capacity.
\Cref{eq:9} establishes that the outbound flow after visiting a GAP equals the inbound flow plus the waste collected from that GAP and, thus, also forbids subtours.

\subsection{Valid inequalities}
\label{subsec:vis}

The model presented above for the MSW~\eqref{tag:m1} is still a difficult problem.
In particular, it contains a lot of \emph{symmetric solutions.}
Two solutions are said to be symmetric if they have the same objective function value but different variable assignments.
Consider the following: during a given day, two trucks undertaking the same collection route would have the same cost.
There is no way for the solver to omit one of them.

One way to address this issue is to add Valid Inequalities (VIs) to the model.
A VI is a constraint that reduces the feasible polytope of the problem without removing every optimal solution.
We decided to focus on VIs for the routing part of the problem because the allocation part is \emph{easy} in comparison.
For examples of VIs in the context of vehicle routing problems, we refer the interested reader to \citet{Dror94}.

One thing to remember is that our graph is asymmetric.
Therefore, we do not need to address symmetries in routes with the same GAPs.
We have developed the following valid inequalities to remove as much symmetry from the optimal solutions as possible.

\subsubsection{Empty start}
A vehicle must start its tour unloaded.
This prevents solutions with different \emph{delivery plans} -- when a vehicle finishes its collection tour below full capacity, we can consider another solution where the vehicle starts with any amount less than the difference.
\begin{equation}
\label{eq:11}
v_{0jlt} = 0, \forall{}\ j\in{I},\ l\in{L},\ t\in{T} \\
\end{equation}

\subsubsection{Vehicle ordering}
We impose that a vehicle with index $l$ can only leave the depot if the vehicle with index $l - 1$ has.
In the case where a solution does not use all available vehicles, we can consider swapping an unused vehicle with a used one.
For brevity, we define: $L' = L \setminus \{0\}$, as the set of vehicles minus the first one.
\begin{equation}
\label{eq:12}
\sum_{i\in{I}} x_{0ilt} \leq{}
\sum_{i\in{I}} x_{0ipt}, \forall\ l \in{L'},\ p = l - 1,\ t\in{T}
\end{equation}

\subsubsection{Furthest visit}
We assign the furthest GAP from the depot to the first vehicle.
Because each GAP must be visited at most once a day, so does the furthest.
Because only one vehicle can visit each GAP on a given day, we can forbid others vehicle than the first vehicle (using $L'$ defined above) to visit the furthest GAP.

\begin{equation}
\label{eq:13}
 \csum{i\in{I}, t\in{T}} x_{ijlt} = 0, \forall\ l \in L', j =\argmax_{i \in I} c_{0i}
\end{equation}

\section{A resolution approach based on Benders decomposition} \label{sec:benders}

\Citet{benders1962partitioning} devised a decomposition method for addressing large MIPs that have a characteristic \emph{block diagonal} structure.
In summary, the method starts by decomposing the \emph{original problem} into a \emph{master problem} and a \emph{subproblem}.
The master problem is a relaxation of the original problem used to determine the values of a subset of its variables.
It is formed by retaining the \emph{complicating variables}, and projecting out the other variables and replacing them with an \emph{incumbent}.
The subproblem is formed around the projected variables and a parameterised version of the complicating variables.
By enumerating the extreme points and rays of the subproblem, the algorithm defines the projected costs and the feasibility requirements, respectively, of the complicating variables.
Because this enumeration is seldom tractable, the algorithm proceeds in the following manner:
\begin{enumerate}
    \item It solves the (relaxed) master problem to optimality, which yields a \emph{candidate solution.}
    \item This candidate solution is used as a parameter in the subproblem.
    \item The resulting problem is solved to optimality and, using LP duality, a set of coefficients are retrieved.
    \item These coefficients are used to generate a constraint, called a ``Benders cut,'' which is added to the master problem.
    \item If the objective function value of the subproblem is equal to the incumbent value in the master problem, the algorithm stops.
    Otherwise, it repeats from point 1.\ using the master problem with the additional constraint.
\end{enumerate}
One key limitation of the classic Benders decomposition is that the subproblem cannot contain integer variables.
This is because of point 3 above: the method needs to use LP duality, which is not well-defined for MIPs.
We use a recent framework called Unified branch-and-Benders-cut \citep[UB\&BC,][]{maheo2020unified} to bypass this issue.
This new framework operates by using a modified B\&C where, at each integer node, it:
\begin{enumerate}
\item solves the LP relaxation of the subproblem to get a lower bound and generate Benders cuts; and,
\item uses a heuristic to determine if the master solution is feasible and, if yes, a valid, global upper bound.
\end{enumerate}
The second point is key: by maintaining a valid upper bound, the framework ensures that no optimal solution is removed during the search.
However, this leads to having a set of \emph{open solutions} after finishing the B\&C tree -- solutions whose objective function value falls between the lower and upper bound.
Thus, the UB\&BC finishes by a \emph{post-processing phase} during which subproblems associated with open solutions are solved to integer optimality.
The combination of maintaining a global upper bound and using a post-processing phase enables the framework to find an optimal solution.

As stated in \cref{sec:mathematical_model}, the problem addressed in this work comprises two characteristic decision-making problems in MSW.
On the one hand, the allocation of bins in the GAPs and, on the other, the design and schedule of routes for the collection vehicles.
This division can be exploited by applying Benders decomposition.
The bins allocation equations are moved to the subproblem while the master problem takes care of designing the schedule and routes of the collection vehicles.

\subsection{Creating the subproblem}

The subproblem allocates bins to each GAP; it is an integer programming problem:%
\named[q(\overline{m}) = \min]{modr2}{SB}
      {\sum_{\substack{i\in{I} \\ u\in{U}}} n_{ui}\ {cin}_{u}}
      {\label{modr2:eq:31}
\sum_{u\in{U}} n_{ui}\ cap_u
& \geq b_{i} \sum_{r\in{R}} \overline{m_{ir}}\ \beta_r & \forall{}\ i\in{I}\\
\label{modr2:eq:32} && \sum_{u\in{U}} n_{ui} & = 1 & \forall\ i\in{I}}
      {n \in \mathbb{B}}
We define the positive continuous variables $\delta_i$ and unrestricted continuous variables $\gamma_i$ as the dual variables of \cref{modr2:eq:31,modr2:eq:32} respectively.
The dual formulation of the LP relaxation of~\eqref{tag:modr2}, which will be used to generate cuts, is then:
\named[q^{LP}(\overline{m}) = \max]{modlp}{LP}
      {\sum_{i \in I} \left(\gamma_i - \delta_i b_i \sum_{r\in R} (\overline{m_{ir}} \beta_r)\right)}
      {\gamma_i - \delta_i \sum_{u \in U} cap_u & \leq \sum_{u \in U} n_{ui} & \forall i \in I}
      {\delta, \gamma \geq 0}

\subsubsection{Heuristic for the subproblem}

In order to apply Benders decomposition when the subproblem has integer variables, an efficient method for solving the subproblem is required.
Therefore, we devised a rounding heuristic procedure based on the LP relaxation of the subproblem:
\begin{enumerate}
\item We solve the LP relaxation of~\eqref{tag:modr2}.
  The (relaxed) solution will contain $n_{ui}$ with fractional values.
\item We estimate the \emph{joint fractional capacity} $K^{f}_{i}$ of each GAP using:
\begin{align}
\label{eq:fractional_capacity}
K^{f}_{i} = \sum_{u\in{U}} n_{ui}\ cap_{u}
\end{align}
\item We define a feasible (non-fractional) bin combination $u \in{U}$ for each GAP by finding the bin combination with minimal cost among those with storage capacity larger than $K^{f}_{i}$. It is guaranteed that there will always be a bin combination which respects this rule since considering \cref{eq:fractional_capacity,modr2:eq:32} implies that: $K^{f}_{i} \leq cap_{u*},\forall\ i\in{I}$, where $u* =\argmax_{u \in{U}} \{cap_{u}\}$.
\end{enumerate}

\subsection{Stating the master problem}

The master problem retains the same constraint structure as~\eqref{tag:m1} but the bin allocation part is replaced by an incumbent variable $q$.
Let us consider the set of extreme points ($\mathcal{O}$) and extreme rays ($\mathcal{F}$) of the LP relaxation of~\eqref{tag:modr2}.
These generate the optimality~\eqref{eq:opticut} and feasibility~\eqref{eq:feascut} cuts, respectively.
Therefore, the master problem is:%
\named{mpb}{MPB}
      {\alpha \csum{i, j \in E(I^0)} \left(c_{ij} + s_i \right) \left(\sum_{\substack{{l\in{L}} \\ t\in{T}}} x_{ijlt} \right) + q}
      {\text{\cref{eq:3,eq:4,eq:5,eq:6,eq:7,eq:8,eq:9,eq:11,eq:12,eq:13}} \span \nonumber \\
&& \sum_{i\in{I}} \left( \gamma_i^f - \delta^f_i b_i \sum_{r\in{R}} (\beta_r m_{ir}) \right) & \leq 0 & \forall\ f \in{\mathcal{F}} \label{eq:feascut} \\
&& \sum_{i\in{I}} \left( \gamma_i^f - \delta^o_i b_i \sum_{r\in{R}} (\beta_r m_{ir}) \right) & \leq q & \forall\ o \in{\mathcal{O}} \label{eq:opticut}}
      {v \geq 0; q \in \mathbb{R}; x \in \mathbb{B}}

\blue{
\subsection{Enhancing the model}

\subsubsection{Partial Benders}

One common issue when using Benders decomposition is that raising the lower bound can take time.
This means that the master problem lacks information on the structure, and thus value, of the subproblem.
This issue has been addressed in stochastic programming with \emph{partial Benders} \citep{Crainic2014}.

The idea is to retain a part of the subproblem's information in the master problem.
This takes the form of adding a relaxed copy of the subproblem's variables.
In our case, we can save a relaxed version of the bin combinations: \cref{modr2:eq:31,modr2:eq:32} with $n \geq 0$.

\subsubsection{Removing symmetric solutions}

Although we use a set of VIs, there is still symmetries in our problem.
Specifically, during the search, we can have a given \emph{master solution} -- in the sense of fixed $m$ variables -- without a corresponding route.
This leads to exploring many solutions with the same bin configuration, or at least having a search that needs to fix a (large) number of routing variables.

To counteract this effect, we would like to have a form of \emph{no-good cut.}
This idea comes from Constraint Programming and consists in forbidding a variable assignment from appearing again.
However, this is not practical in LP.

Instead, we turned towards earlier works in solving integer problems using Benders decomposition.
An \emph{integer L-shaped cut} \citep{Laporte1993} is a constraint that attempts to remove a set of identical master solutions.
It does so by computing the \emph{opportunity cost} of activating a variable.
This cost is only defined for those variables that are active in the current solution.
Finally, this constraint is only active if the solution has the same set of active variables.%
\footnote{ Making it, \emph{de facto}, a no-good cut. }
Do note that this type of constraints only works with binary variables.

The first element we need for our L-shaped cuts is a global lower bound on the problem.
This will allow us to define the opportunity cost as the difference between the subproblem's objective function value and its best possible value.
We define the global lower bound $\mathcal{L}$ as the solution to the subproblem using the least cost bin combination.
That is, we find the smallest possible right-hand side for constraint  \eqref{modr2:eq:31}.
This bin combination can easily be determined by solving:

\named[\mathcal{L} = \min]{min-cost}{LB}
    {\sum_{i \in I} \sum_{r \in R} b_i m_{ir} \beta_r}
    {\sum_{r \in R} m_{ir} & = 1 & \forall i \in I}
    {m \in \mathbb{B}}

Let us denote $\overline{m^{o}}$ the solution at (feasible) iteration $o$ of the Benders approach.
And $S^{o} = \set{(i, r)}{\overline{m_{ir}^{o}} = 1}$ the active variables in a solution .
Using our formulation~\eqref{tag:mpb}, for every feasible subproblem solution $o \in \mathcal{O}$ we have:
\begin{equation}
  \label{eq:l-shaped}
  (q(\overline{m}^{o}) - \mathcal{L}) \left(\sum_{i, r \in S^{o}} m_{ir} - \sum_{i, r \notin S^{o}} m_{ir} \right) - (q(\overline{m}^{o}) - \mathcal{ L }) (|S^{o}| - 1) + \mathcal{ L } \leq q
\end{equation}
}

\section{A working example}
\label{sec:example}

In this section we will present a short example to illustrate our solution approach.
\green{
We will use a toy instance shown in \cref{fig:example-gaps}. The are depicted in the picture correspond to an area of the University neighbourhood of Bahía Blanca retrieved from OpenStreetMap\footnote{https://www.openstreetmap.org/} and, thus, the background image follows the visual code of this engine -- \eg, the red circles represent health facilities. The GAPs' locations are indicated in green circles.
}
It comprises: a two-day-horizon, two GAPs, and two vehicles. We also consider two types of waste bins with a storage capacities of \SI{1.1}{m^3} and \SI{1.73}{m^3} respectively.
\green{These correspond to the two types available in our real-world example.
The costs are derived directly from model \eqref{tag:m1}.}

We will now show the iterations the solution algorithm takes.\footnote{%
    We use the complete problem~\eqref{tag:m1} augmented with valid inequalities \cref{eq:11,eq:12,eq:13}.
    The Benders cuts we generate are ``optimality cuts'' given by~\cref{eq:opticut}.}
At each iteration, we will report:
\begin{itemize}
\item the master solution (routing cost, which includes GAP allocation per vehicle per day);
\item the objective function value of the LP relaxation of the subproblem (lower bound);
\item the heuristic value (upper bound); and,
\item the total cost of the solution.
\end{itemize}
A graphical representation of the iterations is provided in \cref{fig:example-itrs-1,fig:example-itrs-2,fig:example-itrs-3}.
\begin{description}
\item[Iteration 1.]
The first solution uses one vehicles on two days and one vehicle the second day only.
\begin{align*}
  v_{0,0} :\ & (0, 2) \rightarrow (2, 0) \\
  v_{0,1} :\ & (0, 2) \rightarrow (2, 0) \\
  v_{1,0} :\ & (0, 1) \rightarrow (1, 0)
\end{align*}

The routing cost is: $494.2$.
The LP relaxation has an objective function value of $7.96$ while the heuristic has a value of $10.48$.
We add a Benders cut to the master problem and continue.
\item[Iteration 2.]
The second solution uses two vehicles with different routes during one day:
\begin{align*}
  v_0 :\ & (0, 2) \rightarrow (2, 0) \\
  v_1 :\ & (0, 1) \rightarrow (1, 0)
\end{align*}
The routing cost is: $381.8$.
The LP relaxation has an objective function value of $7.06$ while the heuristic has a value of $10.48$.
We add a Benders cut to the master problem and continue.
\item[Iteration 3.]
The third solution found uses a single vehicle with the same route on both days, given by:
\begin{align*}
v_0 :\ &(0, 1) \rightarrow (1, 2) \rightarrow (2, 0)
\end{align*}
The routing cost is: $316.8$.
The LP relaxation has an objective function value of $8.58$ while the heuristic has a value of $10.48$.
We add a Benders cut to the master problem and continue.

At this point, the B\&C will finish as no improving solution can be found, we can progress to the post-processing.
\begin{figure}[ht]
  \centering
\end{figure}
\begin{figure}[ht]
  \centering
  \subfloat[Location of the depot and the two GAPs (green circles) on the toy instance.]{
   \label{fig:example-gaps}
   \includegraphics[width=.46\textwidth]{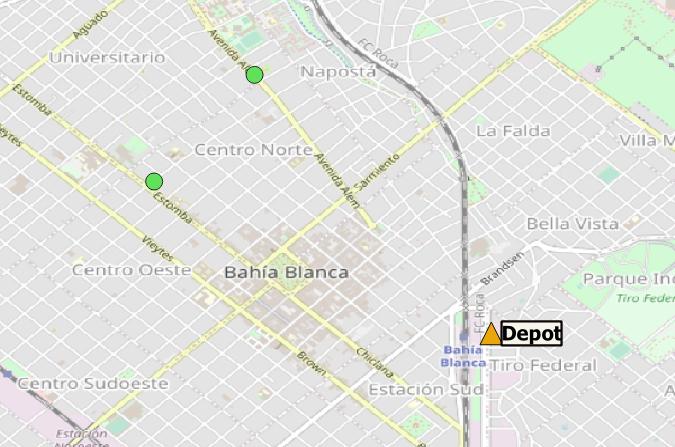}}
  \quad
  \subfloat[It. 1: The first (blue) vehicle uses its route both days, while the second (red) vehicle only operates on the first day.]{
    \label{fig:example-itrs-1}
    \includegraphics[width=.46\textwidth]{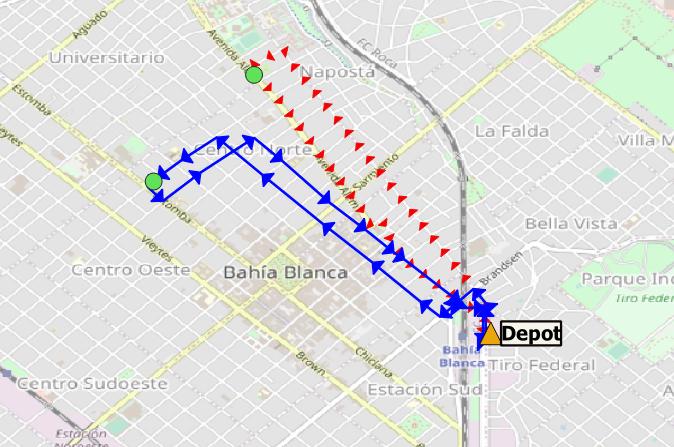}}
  \\
  \subfloat[It. 2: Both vehicles operate during the first day.]{
    \label{fig:example-itrs-2}
    \includegraphics[width=.45\textwidth]{working_example/working_example_b_and_c.jpeg}}
  \quad
  \subfloat[It. 3: Only one vehicle operates during one day]{
    \label{fig:example-itrs-3}
    \includegraphics[width=.46\textwidth]{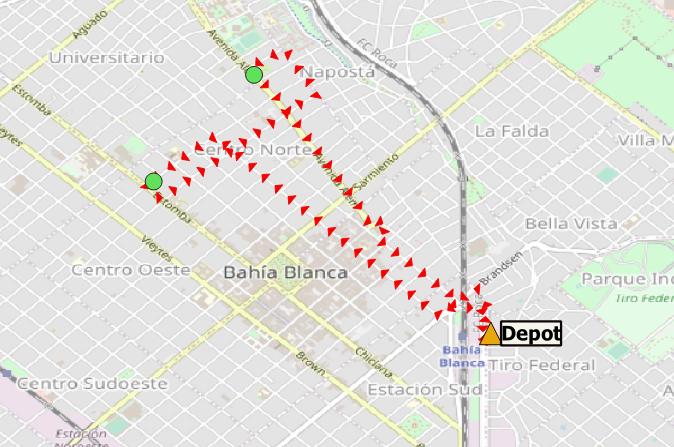}}
    \caption{Working example of the resolution approach.\label{fig:working_example}}
\end{figure}
\item[ Post-processing. ]
At the start of the post-processing phase, the UB\&BC orders solutions according to their lower bound values.
In this case, it will process the solutions in reverse order: 3, 2, 1.

Starting with the solution found in Iteration 3, we solve the subproblem to integer optimality. This gives an optimal value of $10.48$ -- the same as the heuristic.
Being an integer value, it can be used to update the upper bound to: $347.28$ (routing + delivery).

Then, the framework verifies that the remaining open solutions' values are lower that the new-found upper bound.
Both solutions found at Iterations 1 and 2 exceed the best upper bound and are thus skipped.
\end{description}
Our approach has managed to find the optimal solution to the problem.
It did so as an integrated algorithm which solved the routing and allocation problems at once.

\section{Computational experiments}\label{sec:experimentation}

\blue{In this section, we present the computational experimentation of the problem. This includes the description of the real-world instances that were used and implementation details.}

\subsection{Real-world instances}\label{sec:instances}

\begin{figure}
    \centering
    \includegraphics[width=.9\textwidth]{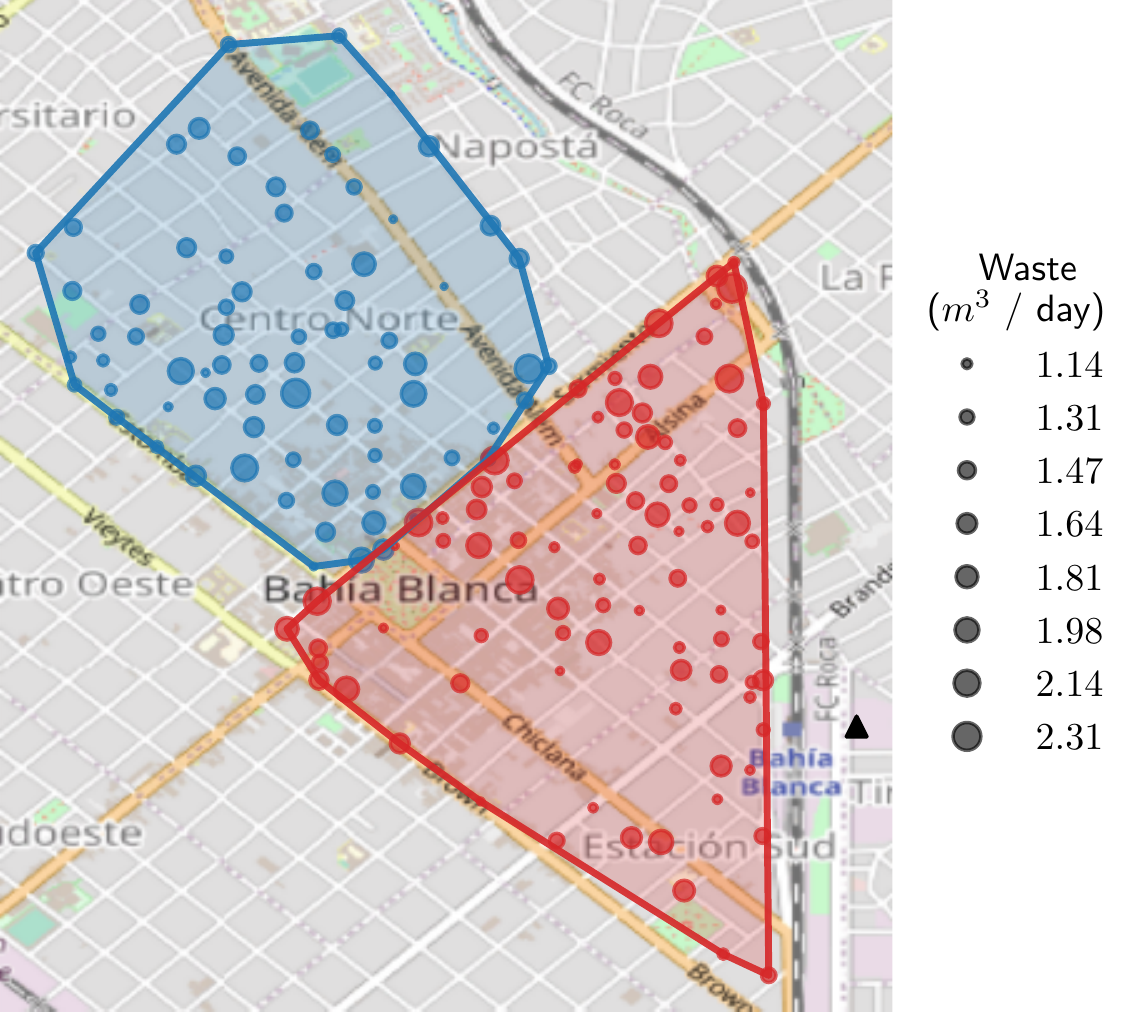}
    \caption{\label{fig:instances}%
    \blue{
        Visualisation of the two instance sets: University (in blue) and Downtown (in red).
        The triangle is the depot.
        The dot size represent the (scaled) daily demand in \si{m^3} per day.}}
\end{figure}

\blue{
Our instances are based on simulated scenarios of the city of Bahía Blanca, Argentina.
Although the city still has a door-to-door collection system, the local government and citizens are interested in more efficient collection systems that allow reductions of the high logistic costs.
For example, using community bins will simplify the logistics costs~\citep{bonomoX2012method,cavallinX2020application,rossitX2020exact}.

The relevant costs of the problem are the bin combination daily cost and the vehicle cost per minute.
Currently, in the door-to-door collection system, the company that performs waste collection in Bah\'{i}a Blanca uses a fleet of rear-loading trucks.
\green{While only small bins can be emptied by a rear-loading truck, almost every bin can be emptied by a side-loading truck.\footnote{To the best of our knowledge, in Argentina the largest bin that allows rear-loading is about \SI{1.1}{m^3}.} In this study we consider three bins that can be emptied by side-loading trucks to allow larger capacities at the GAPs and, thus, taking advantage of using economically convenient bin combinations (\Cref{sec:preproc}).}

We consider three types of commercial side-loading waste bins available in Argentina.
The details are presented in~\cref{tab:bin_types} and were retrieved from surveys to different specialised companies in Argentina.
The life expectancy of bins was estimated to ten years, in line with other similar studies~\citep{brogaard2012quantifying,donzaZ2016full}.
Additionally, the maintenance cost of each bin was estimated at 5\% of the purchasing cost~\citep{donzaZ2016full}.
We calculated the estimated daily cost as the sum of the purchasing and maintenance costs divided by the total amount of days of the expected lifetime.

\begin{table}[]
\small
\centering
\begin{tabular}{l r r r r }
\toprule 
\multirow{2}{*}{Type} & \multirow{2}{*}{\shortstack[c]{Purchase \\ cost (US\$)}} & \multirow{2}{*}{Capacity (\si{m^3})} & \multirow{2}{*}{\shortstack[c]{Occupied \\ area (\si{m^2})}} & \multirow{2}{*}{\shortstack[c]{Estimated daily \\ cost (US\$)}} \\
& & & & \\ \midrule
I & 386.80 & 1.1 & 1.42 & $11.13\exp{-2}$ \\
II & 1102.79 & 2.4 & 2.23 & $31.72\exp{-2}$  \\
III & 1287.24 & 3.2 & 2.60 & $37.03\exp{-2}$  \\
  \bottomrule
\end{tabular}
\caption{\label{tab:bin_types} Details of bin types considered.}
\end{table}

Considering an available space of \SI{5}{m^2} at each GAP, we found eight bin combinations according to the procedure presented in~\cref{tab:bin_combinations}.

\begin{table}[]
\small
\centering
\begin{tabular}{lrrr}
\toprule 
id & Estimated daily cost (US\$) & Capacity (\si{m^3}) & Occupied area (\si{m^2}) \\ \midrule
0       & 0.11  & 1.10  & 1.42  \\
1       & 0.22   & 2.20 & 2.84 \\
2       & 0.32   & 2.40 & 2.23 \\
3       & 0.33   & 3.30 & \green{4.26} \\
4       & 0.43   & 3.50 & \green{3.65} \\
5       & 0.48   & 4.30 & 4.02 \\
6       & 0.64  & 4.80 & 4.46 \\
7       & 0.69  & 5.60 & 4.83 \\
  \bottomrule
\end{tabular}
\caption{\label{tab:bin_combinations} Convenient bin combinations used for the instances.}
\end{table}

We consider an homogeneous fleet of vehicles.
The estimated cost per minute of use ($\alpha$) is taken from the field work of~\citet{donzaZ2016full}.
They estimated that a standard garbage truck with side loader and compactor costs \SI{57.64e-2}[US\$]{} per minute.\footnote{Converted using the official exchange rate of Argentina~\citep{bcra}.}

The service time to empty a GAP -- parameter $s$ of model~\cref{subsec:formulation} -- is not usually considered in the related literature.
However, this time can have a significant impact on the duration of the routes~\citep{giel2021estimating}, especially when there are many GAPs to visit.
In this work, we estimate the time spent collecting waste at a GAP based on the field study performed by~\citet{carlosX2019influence}.
The authors estimated different service time at a GAP depending on factors such as: the collection systems used, the type of trucks, the type of bins, and whether bins have overflowed.
To simplify, in this article we consider that the service time is independent of the bin combination used.
Thus, we take the average of the service time for each bin combinations based on~\citet{carlosX2019influence}.%
\footnote{Performing certain simplifications, the service time of the bin combinations can be estimated assuming that types of bin I, II, and III that are used in this article correspond to the systems S4, S1, and S2 used in~\citet{carlosX2019influence}, respectively.}
The value we used in the computational experimentation is \SI{1.28}{min} (\SI{76.81}{sec} ).

For GAP location and waste generation, we use two datasets from~\citep{cavallinX2020application} who performed a recent field work performed in the city.
These datasets correspond to scenario \emph{F3DM250} for two relevant urban sectors of the city of Bahía Blanca:
\begin{enumerate}
  \item the University neighbourhood, which has a total of 75 GAPs; and,
  \item Downtown, which has 88 GAPs.
\end{enumerate}
Since our problem is much more complex than the problem addressed in~\citet{cavallinX2020application}, which only considers facility location, we use smaller instances.
We generated the smaller instances by picking a random subset of GAP locations using QGIS Random Selection Tool~\citep{QGISsoftware}.
Information about the travel time between GAPs was estimated with Open Source Routing Machine\footnote{\url{http://project-osrm.org/}} using the approach proposed by~\citet{vazquezbrsutX2018ruteo}.

The company that performs collection in Bahía Blanca has a fleet of collection trucks with a capacity of \SI{21}{m^3}.
Since our instances are smaller than the actual collection zones of the city -- the waste equivalent to around seventy GAPs -- we adjust the capacity and size of the fleet in order to not have a trivial instance in which one vehicle can collect all the waste in one trip.
We consider two collection vehicles with capacity set to $$Q = \ceil*{\sum_{i \in I^0} \frac{b_{i}}{2}}$$ where $b_{i}$ are expressed in \si{m^3}.
The fleet size is set to $$|L| = \ceil*{\frac{|I^0|}{2}}$$ where $|I^0|$ is the number of GAPs of the instance plus the depot.
Similarly, the time limit\footnote{That in reality is about five hours.} ($TL$) is downsized to $$TL = \ceil*{\sum_{i, j \in E(I^0)} \frac{c_{ij}}{2}}$$ where $c_{ij}$ are expressed in minutes.

We present the resulting instances in~\cref{tab:instances_description}.%
\footnote{Instances can be retrieved from \url{https://github.com/diegorossit/Set-of-instances-Mah-o-et-al.-2021---ANOR}}
For each instance we report:
\begin{description}
  \item[$|I|$] the number of GAPs considered,
  \item[$|T|$] the number of days of the time horizon,
  \item[$R$] the possible visit combinations considered,
  \item[$|U|$] the number of bin combinations,
  \item[$|L|$] the size of the fleet,
  \item[$Q$] the capacity of the vehicles, and,
  \item[$TL$] the time limit for the routes.
\end{description}
Visit combinations are expressed with the corresponding days when collection is performed within the time horizon. 
Instances with $|U| = 2$ consider the first two bin combinations of \cref{tab:bin_combinations}.

\begin{table}[]
\small
\centering
\begin{tabular}{l r r r r r r r}
\toprule 
Instance & $|I|$ & $|T|$ & $R$ & $|U|$ & $|L|$ & Q (\si{m^3}) & TL (min) \\
\midrule
\multicolumn{8}{c}{University} \\
\midrule
U/5/2/1 & 5 & 2 & \{1,2\}, \{1\}, \{2\} & 2 & 3 & 4 & 55 \\
U/5/2/2 & 5 & 2 & \{1,2\}, \{1\}, \{2\} & 2 & 3 & 3 & 56 \\
U/5/4/1 & 5 & 4 & \{1,2,3,4\}, \{1,3\}, \{2,4\} & 2 & 3 & 4 & 55 \\
U/5/4/2 & 5 & 4 & \{1,2,3,4\}, \{1,3\}, \{2,4\} & 2 & 3 & 3 & 56 \\
U/5/4/3 & 5 & 4 & \{1,2,3,4\}, \{1,3\}, \{2,4\} & 8 & 3 & 4 & 55 \\
U/5/4/4 & 5 & 4 & \{1,2,3,4\}, \{1,3\}, \{2,4\} & 8 & 3 & 3 & 56 \\
U/6/4/1 & 6 & 4 & \{1,2,3,4\}, \{1,3\}, \{2,4\} & 8 & 4 & 4 & 72 \\
U/6/4/2 & 6 & 4 & \{1,2,3,4\}, \{1,3\}, \{2,4\} & 8 & 4 & 4 & 73 \\
U/6/4/3 & 6 & 4 & \{1,2,3,4\}, \{1,3\}, \{2,4\} & 8 & 4 & 4 & 77 \\
U/6/4/4 & 6 & 4 & \{1,2,3,4\}, \{1,3\}, \{2,4\} & 8 & 4 & 4 & 75 \\
U/7/4/1 & 6 & 4 & \{1,2,3,4\}, \{1,3\}, \{2,4\} & 8 & 4 & 5 & 101 \\
U/7/4/2 & 6 & 4 & \{1,2,3,4\}, \{1,3\}, \{2,4\} & 8 & 4 & 5 & 96 \\
\midrule
\multicolumn{8}{c}{Downtown} \\
\midrule
D/5/2/1 & 5 & 2 & \{1,2\}, \{1\}, \{2\} & 2 & 3 & 4 & 47 \\
D/5/2/2 & 5 & 2 & \{1,2\}, \{1\}, \{2\} & 2 & 3 & 4 & 45 \\
D/5/4/1 & 5 & 4 & \{1,2,3,4\}, \{1,3\}, \{2,4\} & 2 & 3 & 4 & 47 \\
D/5/4/2 & 5 & 4 & \{1,2,3,4\}, \{1,3\}, \{2,4\} & 2 & 3 & 4 & 45 \\
D/5/4/3 & 5 & 4 & \{1,2,3,4\}, \{1,3\}, \{2,4\} & 8 & 3 & 4 & 47 \\
D/5/4/4 & 5 & 4 & \{1,2,3,4\}, \{1,3\}, \{2,4\} & 8 & 3 & 4 & 45 \\
D/6/4/1 & 6 & 4 & \{1,2,3,4\}, \{1,3\}, \{2,4\} & 8 & 4 & 5 & 68 \\
D/6/4/2 & 6 & 4 & \{1,2,3,4\}, \{1,3\}, \{2,4\} & 8 & 4 & 4 & 66 \\
D/6/4/3 & 6 & 4 & \{1,2,3,4\}, \{1,3\}, \{2,4\} & 8 & 4 & 5 & 59 \\
D/6/4/4 & 6 & 4 & \{1,2,3,4\}, \{1,3\}, \{2,4\} & 8 & 4 & 5 & 61 \\
D/7/4/1 & 6 & 4 & \{1,2,3,4\}, \{1,3\}, \{2,4\} & 8 & 4 & 5 & 59 \\
D/7/4/2 & 6 & 4 & \{1,2,3,4\}, \{1,3\}, \{2,4\} & 8 & 4 & 5 & 61 \\
  \bottomrule
\end{tabular}
\caption{Instances description.\label{tab:instances_description}}
\end{table}

}

\subsection{Implementation details and execution platform}
\label{subsec:implementation}

The algorithms are implemented in Python 3.5, and we use a UB\&BC framework called \brandec{} \ver.
The solver used is CPLEX v12.7 in its default configuration, we disable heuristics when running the UB\&BC.
We ran the experiments on a computer with Intel Gold 6148 Skylake CPU@2.4GHz and a 8GB RAM limit.

\subsection{Results}
\label{subsec:test_VI}

In this section we present the results of the computational experimentation. We divide the computational experimentation in two parts. In \cref{subsubsec:test_VI} we deal with small instances in order to assess the value of the proposed valid inequalities in the resolution approach. Then, in \cref{subsubsec:tests_with_larger_U} we explore the performance of the proposed Benders approach in comparison to full MIP when solving more complex instances.

\subsubsection{The value of valid inequalities}
\label{subsubsec:test_VI}

\begin{figure}[ht]
\centering
\includegraphics[width=.9\textwidth]{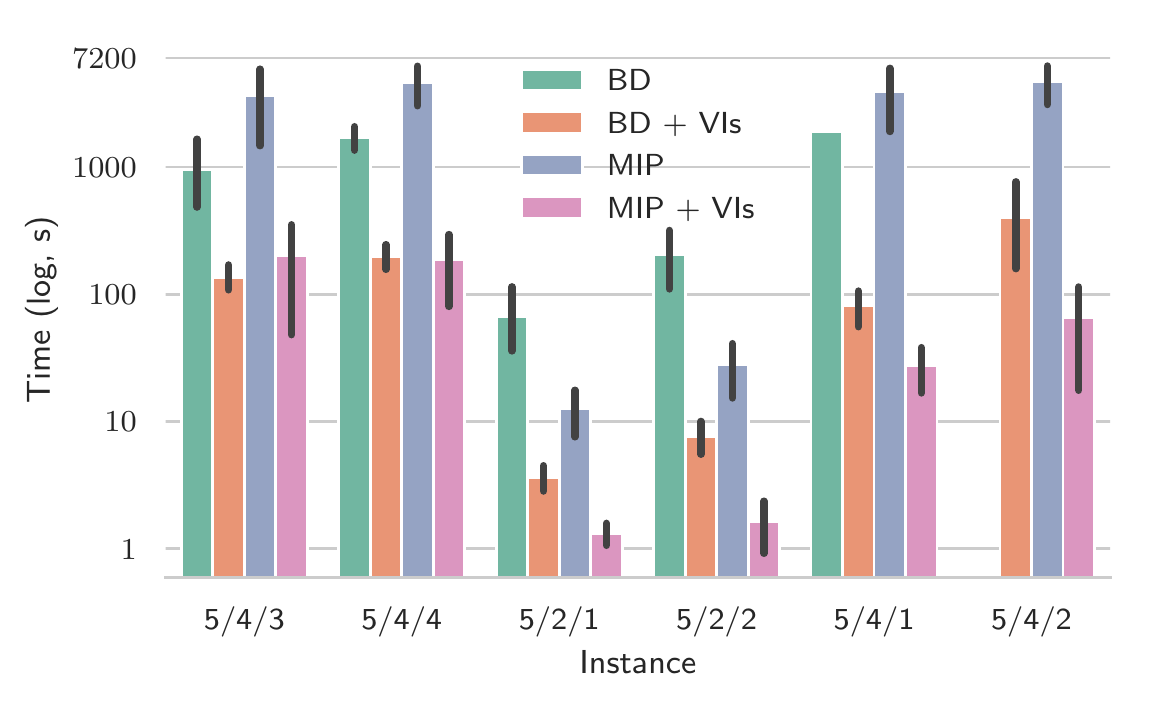}
\caption{\label{fig:vis} Results of using different methods, with or without VIs, to solve a set of reduced instances.}
\end{figure}

\blue{
In order to explore the impact of valid inequalities in the resolution process we solve  the five GAP instances from the both neighbourhoods presented in \cref{tab:instances_description}.
}
\Cref{fig:vis} report the results of solving the resulting problem with:
\begin{description}
\item[MIP:] CPLEX using~\eqref{tag:m1};
\item[MIP + VIs:] CPLEX using~\eqref{tag:m1} augmented with VIs \labelcref{eq:11,eq:12,eq:13};
\item[BD:] our Benders approach; and,
\item[BD + VIs:] our Benders approach augmented with VIs \labelcref{eq:11,eq:12,eq:13}.
\end{description}
We ran five iterations of each configuration and report the minimum solve time.
We can see in \cref{fig:vis} that the VIs are necessary to have reasonable solve times. Both the MIP and our Benders approach benefit from them.
This experiment is not enough to tell for certain whether the Benders approach is better than MIP.

\subsubsection{Reducing symmetry}
\blue{

\green{
We now provide results of using L-shaped cuts in addition to classic Benders cuts.
The main advantage of L-shaped cuts is to reduce symmetry in the master problem by providing a lower bound on solution cost.
In \cref{fig:l-shaped}, we show the runtime when using our Benders approach with or without L-shaped cuts for the set of instances of the University neighbourhood.
}
\begin{figure}[ht]
  \centering
  \includegraphics[width=.9\textwidth]{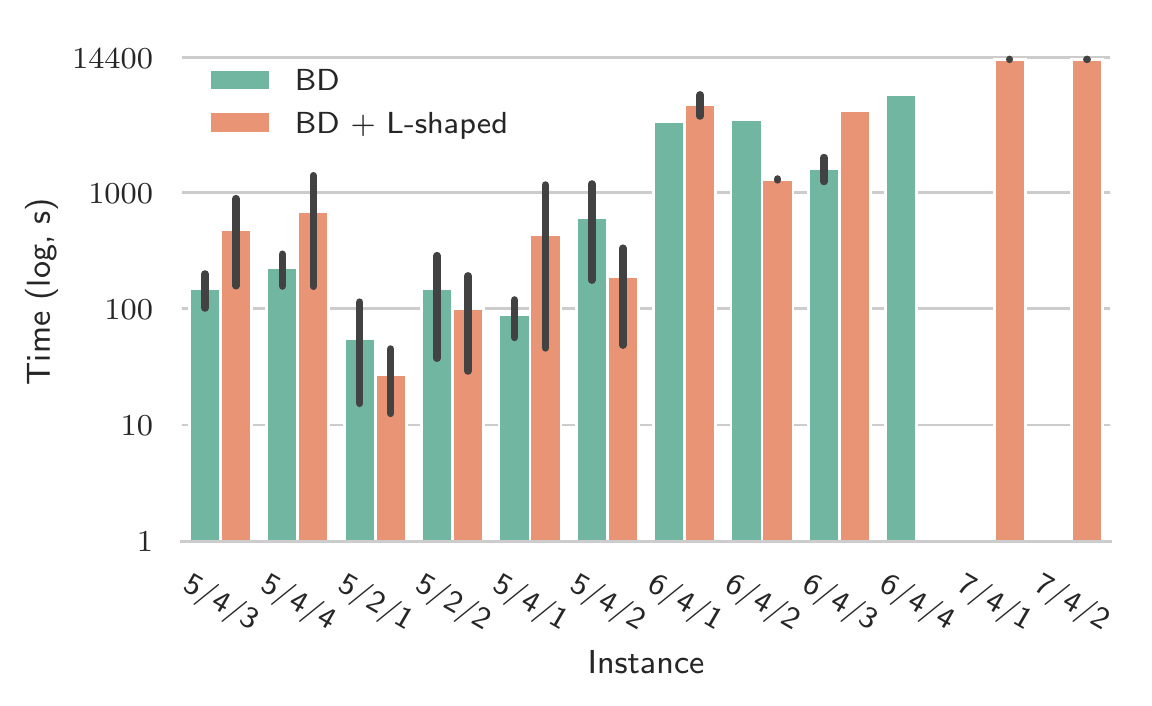}
  \caption{\label{fig:l-shaped}%
    Comparison of our Benders approach with and without L-shaped cuts, we report the total solving time in seconds (log scale).}
\end{figure}

Furthermore, we show the number of master and post-processing iterations, with and without L-shaped cuts, in \cref{tab:l-shaped}.
Overall, L-shaped cuts provide an improvement in the number of iterations required in our Benders approach.
However, this is not a consistent result. 
For master nodes, smaller instances tend to be less affected than larger ones; instances in the university neighbourhood also show less reduction.
The main advantage of L-shaped cuts comes from reducing the number of post-processing iterations. 
By lifting the master solution's objective value early on, we can identify better incumbents.

\begin{table}[ht]
  \small
\begin{tabular}{lrrrrrr}
  \toprule
  & \multicolumn{3}{c}{Master nodes} & \multicolumn{3}{c}{Post-proc.\ nodes} \\
  \cmidrule(lr){2-4} \cmidrule(lr){5-7}
  Instance & Regular & L-Shaped & \% change & Regular & L-Shaped & \% change \\
  \midrule
  \multicolumn{7}{c}{Downtown} \\
5/2/1 & 13 & 47 & 72.3 & 4 & 4 & 0.0 \\
5/2/2 & 41 & 54 & 24.1 & 18 & 11 & -63.6 \\
5/4/1 & 58 & 130 & 55.4 & 40 & 16 & -150.0 \\
5/4/2 & 397 & 94 & -322.3 & 374 & 15 & -2393.3 \\
5/4/3 & 613 & 78 & -685.9 & 608 & 60 & -913.3 \\
5/4/4 & 1795 & 515 & -248.5 & 1792 & 207 & -765.7 \\
6/4/1 & 4437 & 1202 & -269.1 & 288 & 113 & -154.9 \\
6/4/2 & 26260 & 1162 & -2159.9 & 7680 & 678 & -1032.7 \\
6/4/3 & 397 & 1359 & 70.8 & 384 & 181 & -112.2 \\
6/4/4 & 2270 & \textsc{dnf} & n/a & 405 & \textsc{dnf} & n/a \\
7/4/1 & 22 & 16 & -37.5 & 2 & 2 & 0.0 \\
7/4/2 & 22 & 30 & 26.7 & 2 & 3 & 33.3 \\
  \midrule
  \multicolumn{7}{c}{University} \\
5/4/3 & 484 & 468 & -3.4 & 480 & 98 & -389.8 \\
5/4/4 & 482 & 392 & -23.0 & 480 & 121 & -296.7 \\
5/2/1 & 28 & 42 & 33.3 & 4 & 2 & -100.0 \\
5/2/2 & 32 & 239 & 86.6 & 24 & 18 & -33.3 \\
5/4/1 & 84 & 124 & 32.3 & 32 & 5 & -540.0 \\
5/4/2 & 243 & 1508 & 83.9 & 240 & 121 & -98.3 \\
6/4/1 & \textsc{dnf} & 216 & n/a & \textsc{dnf} & 9 & n/a \\
6/4/2 & \textsc{dnf} & \textsc{dnf} & n/a & \textsc{dnf} & \textsc{dnf} & n/a \\
6/4/3 & \textsc{dnf} & \textsc{dnf} & n/a & \textsc{dnf} & \textsc{dnf} & n/a \\
6/4/4 & \textsc{dnf} & \textsc{dnf} & n/a & \textsc{dnf} & \textsc{dnf} & n/a \\
7/4/1 & 4 & 2 & -100.0 & 1 & 1 & 0.0 \\
7/4/2 & 20 & 26 & 23.1 & 3 & 2 & -50.0 \\\bottomrule
\end{tabular}
\caption{\label{tab:l-shaped}%
    Percentage difference in master and post-processing iterations between our Benders approach with and without L-shaped cuts. We report configuration that did not find a solution as \textsc{dnf}.}
\end{table}

\subsubsection{Results over both sectors}
\label{subsubsec:tests_with_larger_U}

We now test our Benders approach against the MIP model on the full set of instances.
\begin{figure}[ht]
  \centering
  \includegraphics[width=.9\textwidth]{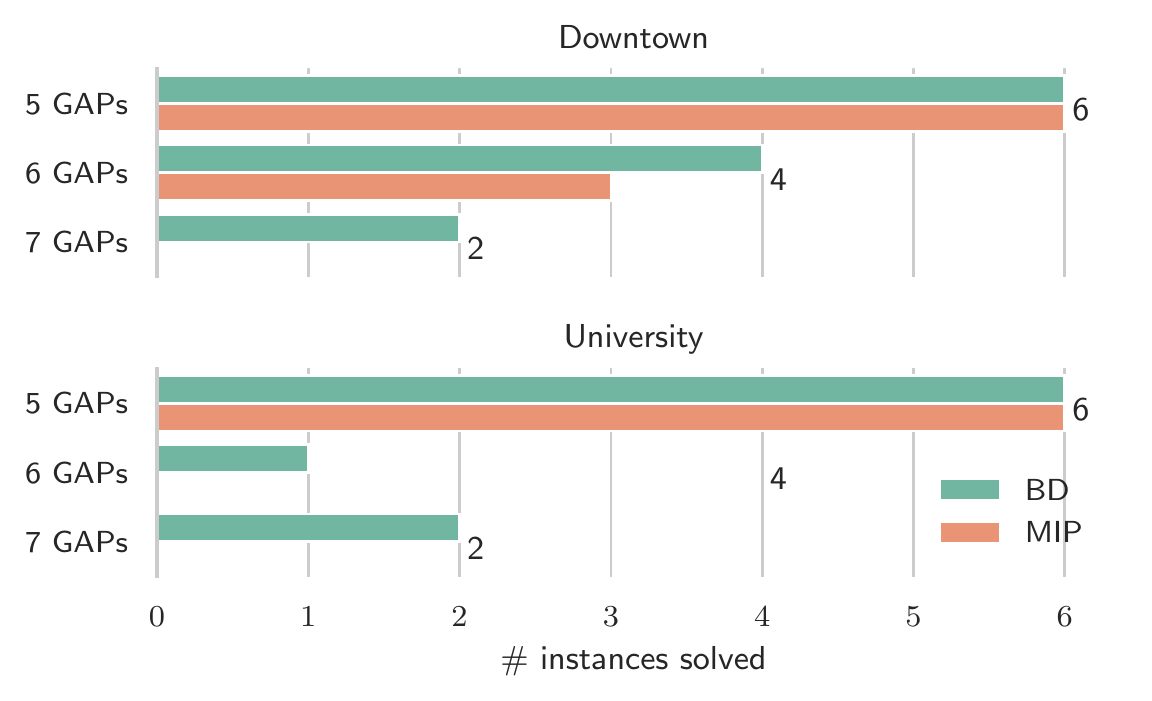}
  \caption{\label{fig:full} Number of instances solved with a MIP or our Benders approach. The black numbers are the total number of instances.}
\end{figure}

In \cref{fig:full} we can see the number of instances solved by each approach.
The MIP approach, solved using CPLEX, does not manage to solve a single instance with 7 GAPs.
With 6 GAPs, the MIP only manages to solve some in the Downtown area.
Again, we can see how the difficulty of instances is not only a function of their size -- GAPs being the most relevant parameter.

}
\section{Conclusion}\label{sec:conc}

\blue{
Municipal solid waste management is a critical issue in modern cities.
Besides the direct environmental and social problems that can arise when it is mishandled, it usually represents a large portion of the municipal budgetary expense~\greenNew{\citep{hoornweg2012waste}}.
Therefore, intelligent decision support tools that can provide high quality of service while also reducing the cost of the system are a major asset for decision makers.

This work addresses two common tactical problems that arise in the reverse logistic chain of solid waste:
\begin{enumerate}
  \item  the design of a pre-collection network, which is based on the location of waste bins; and, 
  \item the routing schedule of collection vehicles, which comprises setting the collection frequency of the bins and the collection routes for the time horizon.
\end{enumerate}
These problems, which are usually solved individually in the related literature, are interdependent in the sense that the solution to one of the problem affects the other.

In this work, we proposed an integrated approach that solves both problems simultaneously, removing trade-offs found in other approaches.
We provided a new MIP formulation, valid inequalities and a resolution approach based on Benders decomposition, using unified branch-and-Benders-cut.
Additionally, we proposed a preprocessing procedure which generates Pareto-optimal waste bin combinations.
This allows reducing the number of integer variables of the MIP formulation. 

Regarding the Benders resolution process, since the subproblem contains integer variables, we devised a heuristic for solving the bin allocation problem.
We performed computational experiments using our approach on a set of real-world instances of two important neighbourhoods of the Argentinean city of Bahía Blanca.
We first tested small instances to show the competitiveness of valid inequalities.
Then, we tested larger instances to analyse the performance of the inclusion of L-shaped cuts in our Benders approach.
The L-shaped cuts allowed our approach to solve the largest instances that we considered (that were not possible to solve with normal Benders).
Finally, we compared the performance of the Benders and MIP approaches on the whole set of instances showing that the proposed Benders approach was more competitive; it was able to solve a larger number of instances within the same time limit.
}

Future work includes expanding computational experiments with larger real-world instances to test the scalability of the approach. \green{Additionally, in this work we have consider some bins that can only be emptied by side-loading trucks. This scenario would imply a replacement of the fleet of collection vehicles in Bah\'{i}a Blanca. Further computational experimentation can be performed to include only bins that can be emptied by rear-loading trucks.}

Another research line is to consider an allocation-first routing-second method.
In that case, the master problem would be comparatively simpler than the subproblem.
Such an approach would require efficient vehicle routing heuristics to work.
We could also explore heterogeneous fleet of vehicles.
Indeed, the city of Bahía Blanca already owns a fleet of vans of small capacity for spot operations.

\begin{acknowledgements}
The second author was supported by the 2018 Australia-Americas PhD Research Internship Program, co-financed by the Australian Academy of Science, the Ministry of Foreign Affairs and Worship of Argentina and the Universidad Nacional del Sur.
\end{acknowledgements}

%
\section*{Conflict of interest}
The authors declare that they have no conflict of interest.

\bibliographystyle{spbasic}      
\bibliography{main_ANOR}   

\end{document}